\documentclass[11pt]{article}
\usepackage[utf8]{inputenc}
\usepackage[T1]{fontenc}
\usepackage[frenchb,english]{babel} 
\usepackage{textcomp}
\usepackage{amsmath,amssymb}
\usepackage{amsthm}
                  
\usepackage{lmodern}

\usepackage[a4paper,margin=2cm]{geometry}

\usepackage{graphicx}             
\usepackage{xcolor}               
\usepackage{microtype,stmaryrd}          

\usepackage{hyperref}
\hypersetup{pdfstartview=XYZ}

\usepackage{bbm}
\usepackage{mathrsfs}
\usepackage{subfig}

\def\build#1_#2^#3{\mathrel{
\mathop{\kern 0pt#1}\limits_{#2}^{#3}}}

\newtheorem{theorem}{Theorem}
\newtheorem{proposition}[theorem]{Proposition}
\newtheorem{definition}[theorem]{Definition}
\newtheorem{lemma}[theorem]{Lemma}

\newtheorem{corollary}[theorem]{Corollary}

\def\w{\mathrm{w}}
\def\t{\mathcal{T}}

\def\ll{\mathcal{L}}
\def\W{\mathcal{W}}
\def\S{\mathcal{S}}

\def\N{\mathbb{N}}

\def\D{\mathbb{D}}
\def\P{\mathbb{P}}
\def\C{\mathbb{C}}
\def\E{\mathbb{E}}
\def\R{\mathbb{R}}
\def\z{\mathcal{Z}}
\def\y{\mathcal{Y}}

\def\ee{\mathcal{E}}
\def\ve{{\varepsilon}}
\def\la{\longrightarrow}
\def\da{\downarrow}

\def\dd{\mathrm{d}}
\def\wh{\widehat}
\def\wt{\widetilde}

\def\bn{\mathbf{n}}
\def\bm{\mathbf{m}}
\def\dg{\mathrm{d}_{\mathrm gr}}

\def\rem{\noindent{\bf Remark. }}

\title{Some explicit distributions for  Brownian motion \\
indexed by the Brownian tree\footnote{Supported by the ERC Advanced Grant 740943 {\sc GeoBrown}}}

\author{Jean-Fran\c cois Le Gall, Armand Riera}

\date{\small Universit\'e Paris-Sud}
\begin{document}
\maketitle

\begin{abstract}
We derive several explicit distributions of functionals of Brownian motion indexed by the Brownian tree. In particular, we give a 
direct proof of a result of Bousquet-M\'elou and Janson identifying the distribution of the density at $0$
of the integrated super-Brownian excursion. 
\end{abstract}
%
%

\section{Introduction}
\label{intro}

The main purpose of the present work is to derive certain explicit distributions for the random process which
we call Brownian motion indexed by the Brownian tree, which has appeared in a variety of different contexts.
As a key tool for the derivation of our main results we use the excursion theory developed in \cite{ALG} for 
Brownian motion indexed by the Brownian tree. In many respects, this excursion theory is similar to the classical 
It\^o theory, which applies in particular to linear Brownian motion and has proved a powerful tool for the calculation
of exact distributions of Brownian functionals.

Let us briefly describe the objects of interest in this work. We define the Brownian tree $\t_\zeta$ as the random compact  $\R$-tree 
coded by a Brownian excursion $\zeta=(\zeta_s)_{s\geq 0}$ distributed according to the (infinite) It\^o measure
of positive excursions of linear Brownian motion. If $\sigma$ stands for the duration of the excursion $\zeta$, this coding means that $\t_\zeta$ is the quotient space of
$[0,\sigma]$ for the equivalence relation defined by $s\sim s'$ if and only if $\zeta_s=\zeta_{s'}=m_\zeta(s,s')$, where 
$m_\zeta(s,s'):=\min\{\zeta_r:s\wedge s'\leq r\leq s\vee s'\}$, and this quotient space is equipped with
the metric induced by $\dd_\zeta(s,s')=\zeta_s+\zeta_{s'}-2m_\zeta(s,s')$. The volume measure $\mathrm{vol}(\dd a)$ on $\t_\zeta$ is defined as 
the push forward of Lebesgue measure on $[0,\sigma]$ under the canonical projection, and the root $\rho$ of $\t_\zeta$ is the
equivalence class of $0$. We note that under the conditioning by $\sigma=1$
(equivalently the total volume is equal to $1$) the tree $\t_\zeta$ is Aldous' Brownian Continuum Random Tree (also called the CRT, see \cite{Al0,Al1}), up to an unimportant scaling factor
$2$. 

Let us turn to Brownian motion indexed by $\t_\zeta$. Informally, given $\t_\zeta$, this is the 
centered Gaussian process $(V_a)_{a\in\t_\zeta}$ such that $V_\rho=0$ and $\mathrm{Var}(V_a-V_b)=\dd_\zeta(a,b)$ for every $a,b\in\t_\zeta$.
This definition is a bit informal since we are dealing with a random process indexed by a {\it random} set. These difficulties can be
overcome easily by using the Brownian snake approach. We let $(W_s)_{s\geq 0}$ be the Brownian snake (whose spatial motion 
is linear Brownian motion started at $0$) driven by the Brownian excursion $(\zeta_s)_{s\geq 0}$. Then, for every $s\geq 0$, 
$W_s$ is a finite path started at $0$ and with lifetime $\zeta_s$, and for every $a\in\t_\zeta$ we may define $V_a$ as the terminal point $\wh W_s$ of the path
$W_s$, for any $s\in[0,\sigma]$ such that $a$ is the equivalence class of $s$ in $\t_\zeta$. The Brownian snake approach thus reduces the
study of a tree-indexed Brownian motion to that of a process indexed by the positive half-line, and we systematically use this approach in the 
next sections.

The total occupation measure $\Theta(\dd x)$ of $(V_a)_{a\in\t_\zeta}$ is the push forward of $\mathrm{vol}(\dd a)$ under the mapping $a\mapsto V_a$,
or equivalently the push forward of Lebesgue measure on $[0,\sigma]$ under $s\mapsto\wh W_s$.
Under the special conditioning $\sigma=1$, this random measure is known as ISE for Integrated Super-Brownian Excursion \cite{Al2} (note that our normalization 
is different from the one in \cite{Al2}).

At this point, we observe that both the Brownian tree (often under special conditionings) and Brownian motion indexed by the Brownian tree
have appeared in different areas of probability theory. The Brownian snake is very closely related to the measure-valued process called 
super-Brownian motion and has proved an efficient tool to study this process (see \cite{Zurich} and the references therein). Super-Brownian motion and ISE arise in
a number of limit theorems for discrete probability models, but also in the theory of interacting particle systems \cite{BCL,CDP,DP} and in a variety of
models of statistical physics \cite{DS,HaS,HoS}. More recently, Brownian motion indexed by the Brownian tree has served as the essential 
building block in the  construction of the universal model of
random geometry called the Brownian map (see in particular  \cite{Uniqueness,Disks,LGM,Mie}). In this connection, we note that the distribution of certain functionals 
of Brownian motion indexed by the Brownian tree is investigated in the article \cite{Del}, which was already motivated by 
asymptotics for random planar maps. 

Let us now explain our main results more in detail. In agreement with the usual notation 
for the Brownian snake, we write $\N_0$ for the (infinite) measure under which 
$(\zeta_s)_{s\geq 0}$ and $(V_a)_{a\in\t_\zeta}$ are defined in the way we just explained --- see Section \ref{preli} for more
details. We are primarily interested in local times, which are the densities of the random measure 
$\Theta(\dd x)$. It follows from the work of Bousquet-M\'elou and Janson \cite{BM0,BMJ}  that $\Theta(\dd x)$ has a continuous density $(\ll^x)_{x\in\R}$
with respect to Lebesgue measure on $\R$, $\N_0$ a.e.
(this fact could also be derived from the earlier work 
of Sugitani \cite{Sug} dealing with super-Brownian motion, see in particular the introduction of \cite{MP}). We also consider the quantity $\sigma_+=\Theta([0,\infty))$ (resp. $\sigma_-=\Theta((-\infty,0])$) corresponding to the
volume of the set of all points $a\in\t_\zeta$ such that $V_a\geq 0$ (resp. $V_a\leq 0$). One of our main technical results (Proposition \ref{distriple}) identifies 
the joint Laplace transform 
$$\N_0(1-\exp(-\lambda \ll^0-\mu_1\sigma_+- \mu_2\sigma_-))\;,\quad \lambda,\mu_1,\mu_2>0,$$
as the solution of the equation $h_{\mu_1,\mu_2}(v)=\sqrt{6}\,\lambda$, where, for every $v\geq 0$,
$$h_{\mu_1,\mu_2}(v)=\sqrt{\sqrt{2\mu_1}+v}\, \Big(2v-\sqrt{2\mu_1}\Big) + \sqrt{\sqrt{2\mu_2}+v}\, \Big(2v-\sqrt{2\mu_2}\Big).$$
In the special case $\mu_1=\mu_2$, this equation can be solved explicitly and leads to the formula
\begin{equation}
\label{disLTsigma}
\N_0\Big(1-\exp(-\lambda \ll^0-\mu\sigma)\Big)=\left\{ \begin{array}{ll}
\sqrt{2\mu}\,\cos\Big(\frac{2}{3}\,\arccos\Big( \frac{\sqrt{3}\,\lambda}{2 (2\mu)^{3/4}}\Big)\Big)\quad&\hbox{if } \frac{\sqrt{3}\,\lambda}{2 (2\mu)^{3/4}}\leq 1,\\
\noalign{\smallskip}
\sqrt{2\mu}\,\cosh\Big(\frac{2}{3}\,\mathrm{arcosh}\Big( \frac{\sqrt{3}\,\lambda}{2 (2\mu)^{3/4}}\Big)\Big)\quad&\hbox{if } \frac{\sqrt{3}\,\lambda}{2 (2\mu)^{3/4}}\geq 1.
\end{array}
\right.
\end{equation}
We can extract the conditional distribution of $\ll^0$ knowing $\sigma$ from the preceding formula. In this way we obtain a short direct proof of a remarkable result of
Bousquet-M\'elou and Janson \cite{BMJ} stating that the local time $\ll^0$ under $\N_0(\cdot\mid \sigma=1)$ (equivalently the
density of ISE at $0$) is distributed as $(2^{3/4}/3)\,T^{-1/2}$, where $T$ is a positive stable
variable with index $2/3$, whose Laplace transform is $\E[\exp(-\lambda T)]=\exp(-\lambda^{2/3})$ (Theorem \ref{lawLT}). The original proof of 
Bousquet-M\'elou and Janson relied on limit theorems for approximations of ISE by discrete labeled trees. Somewhat surprisingly, we are also 
able to obtain an analog of the latter result when instead of conditioning on $\sigma=1$ we condition on $\sigma_+=1$.
Precisely, we get that the local time $\ll^0$ under $\N_0(\cdot\mid \sigma_+=1)$ is distributed as $(2^{9/4}/3)\,D\,T^{-1/2}$, where 
$T$ is as previously and the random variable $D$ is independent of $T$ and 
has density $2x\,\mathbf{1}_{[0,1]}(x)$ (Theorem \ref{lawLT+}). Our proofs are computational and
rely on explicit formulas for moments derived via the Lagrange inversion theorem. It would be interesting to
have more probabilistic proofs and a better understanding of the reason why such simple distributions 
occur.

Because of the connections between the Brownian snake and super-Brownian motion, several of our results
can be restated in terms of distributions of (one-dimensional) super-Brownian motion $(\mathbf{X}_t)_{t\geq 0}$ started from the Dirac measure $\delta_0$.
In particular, we get that the total local time at $0$ (defined as the density at $0$ of the 
measure $\int_0^\infty \dd t\,\mathbf{X}_t$) is distributed as $3^{1/2}2^{-2/3}\,T$ where $T$ is as previously a positive stable variable with index $2/3$
(Corollary \ref{superBM}). 
This is by no means a difficult result (as pointed out to the authors by Edwin Perkins \cite{Per}, the fact that 
the total local time is a stable variable with index $2/3$ can also be derived by a scaling argument, see formula (2.13) in \cite{MP}), but it seems to have remained unnoticed by the specialists of super-Brownian motion. 
The fact that the same variable $T$ occurs in the Bousquet-M\'elou-Janson result suggests the existence of a direct connection between the two 
results, but we have been unable to find such a connection.

The present article is organized as follows. Section \ref{preli} gives a number of preliminaries concerning the Brownian snake.
We have chosen to discuss the Brownian snake with a general spatial motion because it turns out to be useful to
consider also the case where this spatial motion is the pair consisting of a linear Brownian motion and its local time at $0$. 
In fact, Section \ref{secLT} starts with a formula expressing the local time $\ll^0$ in terms of certain exit measures 
of this two-dimensional Brownian snake (Proposition \ref{repreLT}). This expression then leads to an easy calculation of 
the Laplace transform of $\ll^0$, or more generally of $\ll^x$ for any $x\in\R$, under $\N_0$ (Corollary \ref{lawLT-N0}). 
Section \ref{sectriple} gives the key Proposition \ref{distriple} characterizing the joint Laplace transform 
of the triple $(\ll^0,\sigma_+ ,\sigma_-)$ and establishes \eqref{disLTsigma} as a consequence. Finally,
Section \ref{seccondi} derives conditional distributions of the local time $\ll^0$, and also discusses the interpretation
of these distributions in continuous models of random geometry.

\section{Preliminaries}
\label{preli}

\subsection{The Brownian snake}

In this section, we recall some basic facts about the Brownian snake with a general spatial motion. We let 
$\xi$ stand for a Markov process with values in $\R^d$, which starts from $x\in\R^d$ 
under the probability measure $\P_x$. We assume that $\xi$ has continuous sample paths,
and moreover we require the following bound on the increments of $\xi$. There exist three positive constants
$C$, $q>2$ and $\chi>0$ such that for every $t\in[0,1]$ and $x\in \R^d$,
\begin{equation}
\label{moment-condition}
\E_x\Big[\sup_{0\leq s\leq t} |\xi_s - x|^q\Big] \leq C\,t^{2+\chi}.
\end{equation}
Under this moment assumption, we may define the
Brownian snake with spatial motion $\xi$ as a strong Markov process with values in the space
of $d$-dimensional finite paths (see \cite[Section IV.4]{Zurich}). In this work, we will
only need the Brownian snake excursion measures, which we now introduce within
the formalism of snake trajectories \cite{ALG}.

First recall that a ($d$-dimensional) finite path $\w$ is a continuous mapping $\w:[0,\zeta]\la\R^d$, where the
number $\zeta=\zeta_{(\w)}\geq 0$ is called the lifetime of $\w$. We let 
$\W$ denote the space of all finite paths, which is a Polish space when equipped with the
distance
$$d_\W(\w,\w')=|\zeta_{(\w)}-\zeta_{(\w')}|+\sup_{t\geq 0}|\w(t\wedge
\zeta_{(\w)})-\w'(t\wedge\zeta_{(\w')})|.$$
The endpoint or tip of the path $\w$ is denoted by $\wh \w=\w(\zeta_{(\w)})$.
For every $x\in\R^d$, we set $\W_x=\{\w\in\W:\w(0)=x\}$. The trivial element of $\W_x$ 
with zero lifetime is identified with the point $x$ of $\R^d$. 

\begin{definition}
\label{def:snakepaths}
Let $x\in\R^d$. A snake trajectory with initial point $x$ is a continuous mapping $s\mapsto \omega_s$
from $\R_+$ into $\W_x$ 
which satisfies the following two properties:
\begin{enumerate}
\item[\rm(i)] We have $\omega_0=x$ and the number $\sigma(\omega):=\sup\{s\geq 0: \omega_s\not =x\}$,
called the duration of the snake trajectory $\omega$,
is finite (by convention $\sigma(\omega)=0$ if $\omega_s=x$ for every $s\geq 0$). 
\item[\rm(ii)] For every $0\leq s\leq s'$, we have
$\omega_s(t)=\omega_{s'}(t)$ for every $t\in[0,\min_{s\leq r\leq s'} \zeta_{(\omega_r)}]$.
\end{enumerate} 
\end{definition}

If $\omega$ is a snake trajectory, we will write $W_s(\omega)=\omega_s$ and $\zeta_s(\omega)=\zeta_{(\omega_s)}$. 
We denote the set of all snake trajectories with initial point $x$ by $\S_x$.  The set $\S_x$ is equipped with the distance
$$d_{\S_x}(\omega,\omega')= |\sigma(\omega)-\sigma(\omega')|+ \sup_{s\geq 0} \,d_\W(W_s(\omega),W_{s}(\omega'))$$
and the associated Borel $\sigma$-field.

Let $\bn(\dd e)$ denote the classical It\^o measure of positive excursions of linear Brownian motion (see e.g.~\cite[Chapter XII]{RY}). Then $\bn(\dd e)$
is a $\sigma$-finite measure on the space of all continuous functions $s\mapsto e_s$ from $\R_+$ into $\R_+$, and without 
risk of confusion, we will write $\sigma(e)=\sup\{s\geq 0:e_s\not =0\}$, in such a way that we have $0<\sigma(e)<\infty$ and 
$e(s)>0$ for every $0<s<\sigma(e)$, $\bn(\dd e)$ a.e. We consider the usual normalization of $\bn(\dd e)$, so that, 
for every $\ve>0$,
$$\bn\Big(\sup\{e_s:s\geq 0\} >\ve\Big) = \frac{1}{2\ve}.$$
We have then also, for every $\lambda >0$,
\begin{equation}
\label{Laplace-sigma}
\bn(1-\exp(-\lambda \sigma(e))) = \sqrt{\lambda/2},
\end{equation}
and equivalently the distribution of $\sigma(e)$ under $\bn(\dd e)$ is $(2\sqrt{2\pi})^{-1}\,s^{-3/2}\,\dd s$.

\begin{definition}
For every $x\in\R^d$, the Brownian snake excursion measure $\N_x$
is the $\sigma$-finite measure on $\S_x$ characterized by the following two properties:
\begin{enumerate}
\item[\rm(i)] The distribution of $(\zeta_s)_{s\geq 0}$ under $\N_x$ is $\bn$;
\item[\rm(ii)] Under $\N_x$ and conditionally on $(\zeta_s)_{s\geq 0}$, $(W_s)_{s\geq 0}$
is a time-inhomogeneous Markov process whose transition kernels can be described as 
follows: For every $0\leq s\leq s'$,
\begin{enumerate}
\item[$\bullet$] $W_{s'}(t)=W_s(t)$ for all $0\leq t\leq m_\zeta(s,s'):=\min\{\zeta_r:s\leq r\leq s'\}$;
\item[$\bullet$] conditionally on $W_s$, the random path $(W_{s'}(m_\zeta(s,s')+t),0\leq t\leq \zeta_{s'}-m_\zeta(s,s'))$
is distributed as the Markov process $\xi$ started at $W_s(m_\zeta(s,s'))$.
\end{enumerate}
\end{enumerate}
\end{definition}

See again \cite[Chapter IV]{Zurich} for more information about the measures $\N_x$. 
If $F$ is a nonnegative function on $\W_x$, we have the first-moment formula
\begin{equation}
\label{first-mom}
\N_x\Big(\int_0^\sigma F(W_s)\,\dd s\Big) = \E_x\Big[\int_0^\infty F\Big( (\xi_r)_{0\leq r\leq t}\Big) \,\dd t\Big].
\end{equation}

We now turn to exit measures. Let $O$ be an open set in $\R^d$ such that $x\in O$. For every 
$\w\in\W_x$, set
$$\tau_O(\w)=\inf\{t\in[0,\zeta_{(\w)}]: \w(t)\notin O\}$$
with the usual convention $\inf\varnothing =+\infty$. Then $\N_x$ a.e. there exists a random finite measure $\z_O$
supported on $\partial O$
such that, for every bounded continuous function $\varphi$ on $\partial O$, we have
\begin{equation}
\label{approx-exit}
\langle \z_O,\varphi\rangle =\lim_{\ve \to 0} \frac{1}{\ve} \int_0^\sigma \dd s\,\mathbf{1}_{\{\tau_O(W_s)\leq\zeta_s\leq \tau_O(W_s)+\ve\}}\,
\varphi({W}_s(\tau_O(W_s))).
\end{equation}
See \cite[Chapter V]{Zurich}. Then, for every nonnegative measurable function $\varphi$ on $\R^d$,
\begin{equation}
\label{first-mom-exit}
\N_x(\langle \z_O,\varphi\rangle)= \E_x[\varphi(\xi_{\tau_O})\,\mathbf{1}_{\{\tau_O<\infty\}}],
\end{equation}
where in the right-hand side $\tau_O=\inf\{t\geq 0:\xi_t\notin O\}$.

Let us now recall the special Markov property of the Brownian snake, referring to the appendix of \cite{subor} for the proof of a
slightly more precise statement. To this end we consider again the
open set $O$ such that $x\in O$, and fix a snake trajectory $\omega\in\W_x$.  We observe that the set $\{s\geq 0: \tau_O(W_s)<\zeta_s\}$ is open and thus can be 
written as a disjoint union of open intervals $(a_i,b_i)$, $i\in I$ (the indexing set $I$ may be empty if none of the paths $W_s$ exits 
$O$). For every $i\in I$, we may define a new snake trajectory $\omega^{(i)}$ by setting for every $s\geq 0$,
$$\omega^{(i)}_s(t):=\omega_{(a_i+s)\wedge b_i}(\zeta_{a_i}+t)\;,\ \hbox{for every }0\leq t\leq \zeta_{(\omega^{(i)}_s)}:=\zeta_{(a_i+s)\wedge b_i} - \zeta_{a_i}.$$
The snake trajectories $\omega^{(i)}$ represent the excursions of $\omega$ outside $O$ (the word ``outside'' is somewhat misleading since these
excursions typically come back into $O$ though they start on $\partial O$). We also introduce a $\sigma$-field $\mathcal{E}_O$ corresponding informally
to the information given by the paths $W_s$ before they exit $O$ (see \cite{subor} for a more precise definition), and note that $\z_O$ is measurable
with respect to $\mathcal{E}_O$.
Then the special Markov property states that, under $\N_x$ and conditionally on $\mathcal{E}_O$, the point measure
$$\sum_{i\in I} \delta_{\omega^{(i)}}$$
is a Poisson random measure with intensity $\int\z_O(\dd y)\, \N_y(\cdot)$. 

\subsection{Specific properties when $\xi$ is linear Brownian motion}

We finally mention a few more specific properties that hold in the special case where $d=1$ and $\xi$ is 
standard linear Brownian motion. In that case, we have the following scaling property. If $\lambda>0$ and 
\begin{equation}
\label{scaling-rel}
W'_s(t)=\lambda\,W_{s/\lambda^4}(t/\lambda^2)\;,\quad\hbox{for every }0\leq t\leq \zeta'_s:=\lambda^2\zeta_{s/\lambda^4},
\end{equation}
then the distribution of $(W'_s)_{s\geq 0}$ under $\N_x$ is $\lambda^2\N_{\lambda x}$. 

Suppose that the open set $O$ 
is the interval $(-\infty,y)$ with $y>x$, or the interval $(y,\infty)$ with $y<x$. In both cases, the exit measure $\z_O$
can be written as $\z_y\,\delta_{y}$, where $\z_y\geq 0$ and $\delta_y$ denotes the Dirac measure at $y$, and
we have, for every $\lambda>0$,
\begin{equation}
\label{Laplace-exit}
\N_x\Big(1-\exp(-\lambda \z_y)\Big)= \Big(\lambda^{-1/2} + |y-x|\sqrt{2/3}\Big)^{-2}.
\end{equation}
See formula (6) in \cite{CLG}. 

Let $\mathcal{R}:=\{\wh W_s:s\geq 0\}=\{W_s(t):s\geq 0,0\leq t\leq \zeta_s\}$ denote the range of the Brownian snake. Then, for every $y\in \R$, $y\not = x$,
\begin{equation}
\label{range}
\N_x(y\in \mathcal{R})= \frac{3}{2(y-x)^2}= \N_x(\z_y >0).
\end{equation}
See \cite[Section VI.1]{Zurich} for the first equality, and note that the second one follows from \eqref{Laplace-exit}. 

Finally, it follows from the results of \cite{BMJ} that there exists $\N_0$ a.e. a continuous function $(\ll^y)_{y\in \R}$, which is supported on 
$\mathcal{R}$, such that, for every nonnegative measurable function $\varphi$ on $\R$,
$$\int_0^\sigma \dd s \, \varphi(\wh W_s)= \int_\R \dd y\, \varphi(y)\,\ll^y.$$
We call $\ll^y$ the Brownian snake local time at $y$. Note that \cite{BMJ} deals with the case of ISE, that is, with the conditional measure
$\N_0(\cdot\mid \sigma=1)$, but then a scaling argument gives the desired result under $\N_0$.
Next suppose that, for a given $\lambda>0$, $W'$ is defined from $W$ as in \eqref{scaling-rel}. Then, with an obvious
notation, we have $\sigma'=\lambda^4\sigma$ and $\ll'^x=\lambda^3\ll^{x/\lambda}$ for every $x\in\R$, $\N_0$ a.e.
As a consequence, for every $s>0$, the distribution of $\ll^0$ under $\N_0(\cdot \mid \sigma=s)$ is equal to the distribution 
of $s^{3/4}\ll^0$ under $\N_0(\cdot \mid \sigma=1)$.

The scaling property also implies the existence of a constant
$C$ such that, for every $s>0$ and $x\in\R$, we have $\N_0(\ll^x\,|\,\sigma=s)\leq C\,s^{3/4}$ (the case $s=1$ follows
from \cite[Corollary 11.3]{BMJ}, or from a simple argument using Fatou's lemma and the approximation of $\ll^x$ by $(2\ve)^{-1}\int_0^\sigma \dd r\,\mathbf{1}_{\{|\wh W_r-x|<\ve\}}$).

\section{The local time at $0$}
\label{secLT}

In this section and the next ones, we consider the Brownian snake excursion measure $\N_0$ in the case 
where $\xi=B$ is linear Brownian motion. For every $a\in\R$ and $t\geq 0$, we use the notation $L^a_t(B)$ for  the local time of the Brownian motion $B$ at level $a$ and at time $t$. 

For every fixed $s\geq 0$, the path $W_s$ is distributed under $\N_0$ and conditionally on $\zeta_s$
as a linear Brownian path started at $0$ with lifetime $\zeta_s$, and we can define its local time
process at $0$,
$$L^0_t(W_s)=\lim_{\ve\to 0} \frac{1}{\ve}\int_0^t \mathbf{1}_{[0,\ve]}(W_s(r))\,\dd r\;,\quad 0\leq t\leq \zeta_s\,,\ \hbox{a.s.}$$
We may view $L^0(W_s)=(L^0_t(W_s))_{0\leq t\leq \zeta_s}$ as a random element of $\mathcal{W}_0$ with lifetime $\zeta_s$. 
Simple moment estimates show that we can choose a continuous modification of $(L^0(W_s))_{s\geq 0}$ (as a random process
with values in $\mathcal{W}_0$). Moreover, the distribution under $\N_0$ of the two-dimensional process $(W_s,L^0(W_s))_{s\geq 0}$ is the Brownian snake excursion
measure (from the point $(0,0)$ of $\R^2$) for the Markov process $\xi'_t=(B_t,L^0_t(B))$ (note that $\xi'$ may be viewed as a Markov process
with values in $\R^2$, which satisfies \eqref{moment-condition}). 

The point of the preceding discussion is that, under $\N_0$, we can define exit measures for the process $(W_s,L^0(W_s))_{s\geq 0}$
from open subsets $O$ of $\R^2$ containing $(0,0)$, in the way explained in Section \ref{preli}
(e.g. from the approximation formula \eqref{approx-exit}). For every $r>0$, we consider
the exit measure from the open set $O=\R\times(-\infty,r)$ and denote its mass by $\mathscr{X}_r$ (as an application of the
first-moment formula \eqref{first-mom-exit}, this exit measure is a random multiple of the Dirac measure at $(0,r)$).
By convention we also take $\mathscr{X}_0=0$.

 On the other hand, as explained in Section 8 of \cite{ALG}, we can use a famous 
theorem of L\'evy \cite[Theorem VI.2.3]{RY} to give a different presentation of the process $(|W_s|, L^0(W_s))$. To this end, for every $s\geq 0$, write
$$W^\bullet_s(t):= W_s(t)- \min\{W_s(r):0\leq r\leq t\}\,, L^\bullet_s(t):= - \min\{W_s(r):0\leq r\leq t\}\;,\quad\hbox{for }0\leq t\leq \zeta_s.$$
Then the distribution of the pair $(W^\bullet_s,L^\bullet_s)_{s\geq 0}$ under $\N_0$ is equal to the
distribution of  $(|W_s|, L^0(W_s))_{s\geq 0}$ under the same measure. 

Using the preceding identity in distribution of two-dimensional snake trajectories, and 
the approximation \eqref{approx-exit} of exit measures, we get that the process $(\mathscr{X}_r)_{r>0}$
has the same distribution under $\N_0$ as the process $(\z_{-r})_{r>0}$, where we recall that, for every $x\in\R\backslash \{0\}$, $\z_{x}$ denotes the (total mass of the)
exit measure of $(W_s)_{s\geq 0}$ from the open interval $(x,\infty)$ if $x<0$ , or $(-\infty,x)$ if $x>0$ --- of course, by symmetry, $(\z_{-r})_{r>0}$ has the same
distribution as $(\z_{r})_{r>0}$. In particular $\N_0(\mathscr{X}_r>0)=\N_0(\z_r>0)<\infty$ by \eqref{range}. The discussion in \cite[Section 2.4]{LGR} now shows that the process $(\mathscr{X}_r)_{r>0}$
has a c\`adl\`ag modification under $\N_0$, which we consider from now on. Furthermore the distribution of this c\`adl\`ag modification under $\N_0$
can be interpreted as the excursion measure of the continuous-state branching process with branching mechanism $\phi(u)=\sqrt{8/3}\,u^{3/2}$
(the $\phi$-CSBP in short, see \cite[Chapter II]{Zurich} for a brief presentation of continuous-state branching processes). This means that, if $\alpha >0$, and 
$\sum_{i\in I} \delta_{\omega_i}$
is a Poisson point measure with intensity $\alpha\,\N_0$, the process $Y$ defined by  $Y_0=\alpha$ and $$Y_r=\sum_{i\in I} \mathscr{X}_r(\omega_i)$$
for every $r>0$,
 is a $\phi$-CSBP started from $\alpha$. Note that the right-hand side of the last display is a finite sum since $\N_0(\mathscr{X}_r>0)<\infty$. 

Recall our notation $\mathcal{L}^0$ for the Brownian snake local time at $0$. 

\begin{proposition}
\label{repreLT}
We have
$$\mathcal{L}^0=\int_0^\infty \dd r\,\mathscr{X}_r\;,\quad \N_0\hbox{ a.e.}$$
\end{proposition}

This proposition is obviously related to the identity (37) in \cite[Proposition 25]{LGR}, which is however 
concerned with the local time $\ll^x$ at a level $x>0$. Unfortunately, the case $x=0$ seems to require a 
different argument.

\proof It will be convenient to write $\wh L(W_s)=L^0_{\zeta_s}(W_s)$ and
$$L^*=\max\{\wh L(W_s):0\leq s\leq \sigma\}.$$
For every $\ve >0$, set
$$\ll^{0,\ve}:= \ve^{-1}\int_0^\sigma \dd s\,\mathbf{1}_{\{0<\wh W_s<\ve\}}.$$
Then $\ll^{0,\ve}\la\ll^0$ as $\ve\to 0$, $\N_0$ a.e.
We also introduce, for every fixed $\delta>0$,
$$\ll^{0,\ve,(\delta)}:= \ve^{-1}\int_0^\sigma \dd s\,\mathbf{1}_{\{0<\wh W_s<\ve, \wh L(W_s)> \delta\}}.
$$
We observe that, for every $\ve,\delta>0$, we can use the first-moment formula \eqref{first-mom} to compute
\begin{equation}
\label{repreLT-tech1}
\N_0(\ll^{0,\ve}-\ll^{0,\ve,(\delta)})
= \N_0\Big( \ve^{-1}\int_0^\sigma \dd s\,\mathbf{1}_{\{0<\wh W_s<\ve,  \wh L(W_s)\leq \delta\}}\Big)
= \ve^{-1}\E_0\Big[ \int_0^\infty \dd t\,\mathbf{1}_{\{0<B_t<\ve, L^0_t(B)\leq\delta\}}\Big]=\delta,
\end{equation}
where the last equality follows from a standard Ray-Knight theorem for Brownian local times
\cite[Theorem IX.2.3]{RY}. We then need the following lemma.

\begin{lemma}
\label{repreLT-tech}
For every $\delta>0$, we have 
$$\lim_{\ve \to 0} \ll^{0,\ve,(\delta)} = \int_\delta^\infty \dd r\,\mathscr{X}_r\;,$$
in probability under $\N_0(\cdot\mid L^*\geq \delta)$.
\end{lemma}

Let us postpone the proof of this lemma and complete that of Proposition \ref{repreLT}. Write 
$\wt \ll^0=\int_0^\infty \dd r\,\mathscr{X}_r$ and $\wt \ll^{0,(\delta)}=\int_\delta^\infty \dd r\,\mathscr{X}_r$
to simplify notation, and for $a>0$ set $\N_0^{(a)}=\N_0(\cdot\mid L^*\geq a)$. Then, for every $\alpha>0$,
\begin{align}
\label{repreLT-tech2}
\N_0^{(a)}(|\ll^0-\wt \ll^0|>\alpha)
&\leq 
\N_0^{(a)}(|\ll^0- \ll^{0,\ve}|>\alpha/4)
+\N_0^{(a)}(|\ll^{0,\ve}-\ll^{0,\ve,(\delta)}|>\alpha/4)\nonumber\\
&\ + \N_0^{(a)}(|\ll^{0,\ve,(\delta)}-\wt \ll^{0,(\delta)}|>\alpha/4)
+ \N_0^{(a)}(|\wt \ll^{0,(\delta)}-\wt \ll^{0}|>\alpha/4).
\end{align}
Let $\gamma>0$. We can fix $\delta>0$ small enough so that, for every $\ve >0$, the second and the fourth term in the
right-hand side of \eqref{repreLT-tech2} are smaller than $\gamma/4$ (we use \eqref{repreLT-tech1}
for the second term). Then, if $\ve>0$ is small enough, the first and the third term are also smaller than $\gamma/4$
(using Lemma \ref{repreLT-tech} for the third term). We conclude that $\N_0^{(a)}(|\ll^0-\wt \ll^0|>\alpha)\leq \gamma$
and since $\alpha$ and $\gamma$ were arbitrary this gives the desired result $\wt\ll^0=\ll^0$. \endproof

\medskip
\noindent{\it Proof of Lemma \ref{repreLT-tech}.} We keep the notation $\wt \ll^{0,(\delta)}$ introduced in the previous proof.
We first observe that
\begin{equation}
\label{repreLT-1}
\wt \ll^{0,(\delta)}=\lim_{\ve \to 0} \ve \sum_{k=0}^\infty \mathscr{X}_{\delta +k\ve}\;,\ \N_0\hbox{ a.e.}
\end{equation}
and on the other hand,
\begin{equation}
\label{repreLT-2}
\ll^{0,\ve,(\delta)} = \ve^{-1}\sum_{k=0}^\infty\mathcal{H}^{\ve,(\delta)}_k,
\end{equation}
where 
$$\mathcal{H}^{\ve,(\delta)}_k=\int_0^\sigma \dd s\,\,\mathbf{1}_{\{0<\wh W_s<\ve, \;\delta+k\ve<\wh L(W_s)\leq \delta+(k+1)\ve\}}.$$
The idea of the proof is to bound $\N_0(|\ve\mathscr{X}_{\delta +k\ve} - 
\ve^{-1} \mathcal{H}^{\ve,(\delta)}_k|)$, for every fixed $k\geq 0$.
To this end, we apply the special Markov property to the Brownian snake with spatial motion $(B_t,L^0_t(B))$ and the
open set $O=\R\times (-\infty,\delta+k\ve)$, noting that the event $\{L^*\geq \delta\}$ is then $\ee_O$-measurable. It follows that,
under $\N_0(\cdot\mid L^*\geq\delta)$ and
conditionally on $\mathscr{X}_{\delta +k\ve}=a$, the quantity $\mathcal{H}^{\ve,(\delta)}_k$
is distributed as 
$$\int \mathcal{N}(\dd \omega)\,\mathcal{U}_\ve(\omega)$$
where $\mathcal{N}(\dd \omega)$ is a Poisson point measure with intensity $a\,\N_0$, and the random variable $\mathcal{U}_\ve$
is defined under $\N_0$ by
$$\mathcal{U}_\ve = \int_0^\sigma \dd s\,\,\mathbf{1}_{\{0<\wh W_s<\ve,0< \wh L(W_s)\leq \ve\}}.$$
Hence, conditionally on $\mathscr{X}_{\delta +k\ve}=a$, $\mathcal{H}^{\ve,(\delta)}_k$ has the distribution of
$U^\ve_a$, where $(U^\ve_t)_{t\geq 0}$ is the subordinator whose L\'evy measure is the distribution of $\mathcal{U}_\ve$
under $\N_0$. Note that $\E[U^\ve_1]=\N_0(\mathcal{U}_\ve)=\ve^2$ by \eqref{repreLT-tech1}.

By a scaling argument, we get that $(U^\ve_t)_{t\geq 0}$ has the same distribution as $(\ve^4\,U^1_{\ve^{-2}t})_{t\geq 0}$. Next the law
of large numbers shows that
\begin{equation}
\label{LLN}
\lim_{t\to\infty} \sup_{s\leq t} \E\Big[ \frac{|U^1_s -s|}{t}\Big]=0.
\end{equation}
Fix $A>0$ and consider the event $E_A:=\{L^*\leq A\}\cap\{\sup\{\mathscr{X}_{r}:r\geq 0\}\leq A\}$. Notice that on this event we have
$\mathscr{X}_{\delta +k\ve}=0$ and $\mathcal{H}^{\ve,(\delta)}_k=0$ as soon as $\delta+k\ve >A$. It follows that
\begin{align*}
&\N_0\Big(\mathbf{1}_{E_A}\,|\ve \sum_{k=0}^\infty \mathscr{X}_{\delta +k\ve} - \ve^{-1}\sum_{k=0}^\infty\mathcal{H}^{\ve,(\delta)}_k|
\;\Big| \,L^*\geq \delta\Big)\\
&\quad\leq \ve (\lfloor A/\ve\rfloor +1) \sup_{0\leq k\leq \lfloor A/\ve\rfloor}
\N_0\Big(\mathbf{1}_{\{\mathscr{X}_{\delta +k\ve}\leq A\}} \, |\mathscr{X}_{\delta +k\ve} -\ve^{-2}\mathcal{H}^{\ve,(\delta)}_k|\;\Big|\; L^*\geq \delta\Big)\\
&\quad\leq \ve (\lfloor A/\ve\rfloor +1) \sup_{0\leq a\leq A} \E[|\ve^{-2}U^\ve_a -a|]\\
&\quad= \ve (\lfloor A/\ve\rfloor +1) \sup_{0\leq s\leq A/\ve^2} \E[\ve^2|U^1_s-s|],
\end{align*}
which tends to $0$ as $\ve\to 0$, by \eqref{LLN}. The statement of the lemma follows, recalling \eqref{repreLT-1}
and \eqref{repreLT-2}. \hfill$\square$

\begin{corollary}
\label{lawLT-N0}
For every $\lambda>0$,
\begin{equation}
\label{Laplace-local}
\N_0(1-e^{-\lambda\ll_0})=\frac{3^{1/3}}{2}\,\lambda^{2/3}.
\end{equation}
The distribution of $\ll^0$ under $\N_0$ has density
$$h(\ell)= \frac{3^{-2/3}}{\Gamma(1/3)}\,\ell^{-5/3}$$
with respect to Lebesgue measure on $(0,\infty)$.
\end{corollary}

\proof Let $(X_r)_{r\geq 0}$ denote a $\phi$-CSBP started from $1$, where we recall that $\phi(u)=\sqrt{8/3}\,u^{3/2}$. Using the interpretation of the distribution of 
$(\mathscr{X}_r)_{r>0}$ under $\N_0$ and the exponential formula for Poisson measures, we have
$$\E \Big[\exp\Big(-\lambda \int_0^\infty \dd r\,X_r\Big)\Big]= \exp\Big(-\N_0\Big(1-\exp\Big(-\lambda\int_0^\infty \dd r\,\mathscr{X}_r\Big)\Big)\Big).$$
It then follows from Proposition \ref{repreLT} that
$$\N_0(1-e^{-\lambda\ll_0})= -\log \E \Big[\exp\Big(-\lambda \int_0^\infty \dd r\,X_r\Big)\Big].$$
The classical Lamperti transformation
\cite{CLU,Lam0} shows that $\int_0^\infty \dd r\,X_r$ has the same distribution as $T_0:=\inf\{t\geq 0: Y_t=0\}$,
where $(Y_t)_{t\geq 0}$ denotes a stable L\'evy process with no negative jumps started from $1$, whose distribution
is characterized by the Laplace transform $\E[\exp(-\lambda (Y_t-1))]=\exp(t\,\phi(\lambda))$. It is then classical
(see e.g.~\cite[Chapter VII]{Ber}) that
$$\E[e^{-\lambda T_0}]=e^{-\phi^{-1}(\lambda)},$$
where $\phi^{-1}(\lambda)=(3/8)^{1/3}\,\lambda^{2/3}$ is the inverse function of $\phi$. This completes the proof
of the first assertion. The density of $\ll^0$ is then obtained by inverting the Laplace transform. \endproof

In the next corollary, we consider a one-dimensional super-Brownian motion $(\mathbf{X}_t)_{t\geq 0}$
with quadratic branching mechanism $\psi(u)=2u^2$ (the choice of the constant $2$ is only for convenience,
since a scaling argument will give a similar result with a general quadratic branching mechanism). Then it is well known that
we can define the associated (total) local times  as the unique (random) continuous function $(\mathbf{L}^a)_{a\in\R}$ such that
$$\int_0^\infty \dd t\,\langle \mathbf{X}_t,f\rangle = \int_{\R} \dd a\,f(a)\,\mathbf{L}^a,$$
for every Borel function $f:\R\la\R_+$.
See in particular Sugitani \cite{Sug}.

\begin{corollary}
\label{superBM}
Suppose that $\mathbf{X}_0=\alpha\,\delta_0$ for some $\alpha>0$. Then, for every $a\in\R$ and $\lambda >0$,
\begin{equation}
\label{identLT01}
\E[e^{-\lambda \mathbf{L}^a}]=\exp\left(-\alpha \,\frac{3^{1/3}}{2}\Big(\lambda^{-1/3} +3^{-1/3}\,|a|\Big)^{-2}\right).
\end{equation}
In particular,
\begin{equation}
\label{identLT00}
\E[e^{-\lambda \mathbf{L}^0}]= \exp\left(-\alpha\frac{3^{1/3}}{2}\,\lambda^{2/3}\right),
\end{equation}
so that $\mathbf{L}^0$ is a positive stable variable with index $2/3$. 
\end{corollary}

\proof We rely on the Brownian snake construction of super-Brownian motion (see in
particular \cite[Chapter 4]{Zurich}). We may assume that $(\mathbf{X}_t)_{t\geq 0}$ is constructed in
such a way that there exists a Poisson point measure $\mathcal{N}=\sum_{i\in I}\delta_{\omega_i}$
with intensity $\alpha \N_0$, such that, for every Borel function $f:\R\la\R_+$,
$$\int_{\R} \dd a\,f(a)\,\mathbf{L}^a= \int_0^\infty \dd t\,\langle \mathbf{X}_t,f\rangle
=\sum_{i\in I} \int_0^{\sigma(\omega_i)} \dd s\,f(\wh W_s(\omega_i))= \sum_{i\in I} \int_{\R} \dd a\,f(a)\,\ll^a(\omega_i).$$
It follows that we have
\begin{equation}
\label{identLT}
\mathbf{L}^a= \sum_{i\in I} \ll^a(\omega_i)
\end{equation}
for Lebesgue a.e. $a\in\R$. The left-hand side is continuous in $a$, and the right-hand side is continuous on $\R\backslash\{0\}$ since, for every 
$\delta>0$, there are only finitely many $i\in I$ such that $\mathcal{L}^a(\omega_i)$ is nonzero for some $a$ with $|a|>\delta$. So
\eqref{identLT} holds for every $a\in\R\backslash\{0\}$. In fact it is easy to see that \eqref{identLT} also holds for $a=0$. 
First note that, by Fatou's lemma, $\mathbf{L}^0\geq  \sum_{i\in I} \ll^0(\omega_i)$, so that it suffices to check that
$$\E[e^{-\mathbf{L}^0}]= \E\Big[\exp\Big(-\sum_{i\in I} \ll^0(\omega^i)\Big)\Big].$$
The left-hand side is the limit when $a\to 0$ of $\E[e^{-\mathbf{L}^a}]=\exp(-\N_0(1-e^{-\ll^a}))$ and the right-hand side 
is equal to $\exp(-\N_0(1-e^{-\ll^0}))$. So we only need to verify that $\N_0(1-e^{-\ll^a})$ tends to $\N_0(1-e^{-\ll^0})$
as $a\to 0$, which is easy by conditioning on $\sigma$ and then using the bound $\N_0(1-e^{-\ll^a}\,|\,\sigma=s)\leq C(s^{3/4}\wedge 1)$
to justify dominated convergence.

Formula \eqref{identLT00} follows from the case $a=0$ of \eqref{identLT} as an immediate application of \eqref{Laplace-local} and the
exponential formula for Poisson measures. As for formula \eqref{identLT01}, it is enough to verify that
\begin{equation}
\label{identLT3}
\N_0(1-e^{-\lambda\ll^a}) = \frac{3^{1/3}}{2}\Big(\lambda^{-1/3} +3^{-1/3}\,|a|\Big)^{-2}.
\end{equation}
Fix $a>0$ for definiteness, and recall our notation $\z_a$ for the total mass of the exit measure from $(-\infty,a)$. Write $(\omega'_j)_{j\in J}$
for the excursions of the Brownian snake outside $(-\infty,a)$. By the special Markov property, under $\N_0$ and conditionally on $\z_a$, the
point measure $\sum_{j\in J} \delta_{\omega'_j}$ is Poisson with intensity $\z_a\,\N_a$. Moreover, the first part of the 
proof shows that we have $\mathbf{L}^a=\sum_{j\in J} \ll^a(\omega'_j)$,  $\N_0$ a.e., and therefore
$$\N_0(1-e^{-\lambda\ll^a}) = \N_0\Big(1- \exp\Big(-\z_a \N_0(1-\exp(-\lambda \ll^0)\Big)\Big).$$
Then \eqref{identLT3} 
follows from  \eqref{Laplace-local} and \eqref{Laplace-exit}. \endproof

\rem An alternative way to derive the previous two corollaries would be to use the known connections between
super-Brownian motion or the Brownian snake and partial differential equations. See formula (1.13) in \cite{MP}, and note that,
as a function of $a$, the right-hand side of \eqref{identLT3} solves the differential 
equation $\frac{1}{2}u''= 2u^2 -\lambda \delta_0$ in the sense of distributions. On the other hand, our 
method provides a better probabilistic understanding of the results and the derivation 
of \eqref{Laplace-local} in particular relies on Proposition \ref{repreLT} which is of independent interest and will
play a significant role in the proofs of the next section.

\section{The joint distribution of the local time and the time spent above and below $0$}
\label{sectriple}

Our next goal is to discuss the joint distribution of $(\ll^0,\sigma_+,\sigma_-)$ under $\N_0$, where we write
$$\sigma_+:=\int_0^\sigma \mathbf{1}_{\{\wh W_s>0\}}\,\dd s\;,\quad \sigma_-:=\int_0^\sigma \mathbf{1}_{\{\wh W_s<0\}}\,\dd s.$$

\begin{proposition}
\label{distriple}
Let $\lambda,\mu_1,\mu_2\geq 0$, and consider the function $h_{\mu_1,\mu_2}:[0,\infty)\la \R$ defined by
$$h_{\mu_1,\mu_2}(v)=\sqrt{\sqrt{2\mu_1}+v}\, \Big(2v-\sqrt{2\mu_1}\Big) + \sqrt{\sqrt{2\mu_2}+v}\, \Big(2v-\sqrt{2\mu_2}\Big) .$$
Then the quantity 
$$v(\lambda,\mu_1,\mu_2):= \N_0(1-\exp(-\lambda \ll^0-\mu_1\sigma_+-\mu_2\sigma_-))$$
is the unique solution of the equation $h_{\mu_1,\mu_2}(v)=\sqrt{6}\, \lambda$. 
%
\end{proposition}

\proof First note that the quantities $v(\lambda,\mu_1,\mu_2)$ are finite, since $v(\lambda,\mu,\mu)\leq \N_0(1-\exp(-\lambda\ll^0))
+ \N_0(1-\exp(-\mu\sigma))<\infty$ by \eqref{Laplace-sigma} and \eqref{Laplace-local}. Then,
suppose that, under the probability measure $\P$, we are given a sequence $(\eta_i)_{i\geq 0}$ of 
independent Bernoulli variables with parameter $1/2$, and a sequence $(U_i)_{i\geq 0}$
of i.i.d. nonnegative random variables with density $(2\pi u^5)^{-1/2}\exp(-1/2u)$ for $u>0$. 
We note that, for every $\beta>0$, we have
\begin{equation}
\label{LaplaceU}
\E[\exp(-\beta U_1)]= (1+\sqrt{2\beta})\,\exp(-\sqrt{2\beta}).
\end{equation}
The reason for introducing these two sequences is the following fact. If $(t_i)_{i\geq 0}$ is a measurable
enumeration of the jump times of the process $(\mathscr{X}_t)_{t\geq 0}$ (under $\N_0$), the 
conditional distribution of the pair $(\sigma_+,\sigma_-)$ under $\N_0$ and knowing $(\mathscr{X}_t)_{t\geq 0}$
is the law of 
$$\Big( \sum_{i=0}^\infty \eta_i\, U_i\,(\Delta\mathscr{X}_{t_i})^2,  \sum_{i=0}^\infty (1-\eta_i)\, U_i\,(\Delta\mathscr{X}_{t_i})^2\Big).$$
This fact is a consequence of the excursion theory developed in \cite{ALG} (in particular Theorem 4 and Proposition 31 of \cite{ALG}).
In this theory, excursions away from $0$ are in one-to-one correspondence with the jumps of $(\mathscr{X}_t)_{t\geq 0}$, so
that in the preceding display $\eta_i$ gives the sign of the associated excursion ($\eta_i=1$ for a positive excursion and
$\eta_i=0$ for a negative one), and $U_i\,(\Delta\mathscr{X}_{t_i})^2$ corresponds to the duration of this excursion. We refer to
\cite{ALG} for more details. 

Using also Proposition \ref{repreLT} and \eqref{LaplaceU}, it follows that
$$\N_0\Big(\exp(-\lambda\ll^0 -\mu_1 \sigma_+-\mu_2\sigma_-)\,\Big|\,(\mathscr{X}_t)_{t\geq 0}\Big)=
\exp\Big(-\lambda\int_0^\infty \dd t\, \mathscr{X}_t\Big)\,\prod_{i=0}^\infty F(\mu_1,\mu_2, (\Delta\mathscr{X}_{t_i})^2),$$
where we have set, for every $x>0$, 
$$F(\mu_1,\mu_2,x):= \frac{1}{2}\Big( (1+\sqrt{2\mu_1 x})\,\exp(-\sqrt{2\mu_1 x}) + (1+\sqrt{2\mu_2 x})\,\exp(-\sqrt{2\mu_2 x})\Big).$$
Hence, with the notation of the theorem, we have
$$v(\lambda,\mu_1,\mu_2)=\N_0\Big( 1 - \exp\Big(-\lambda\int_0^\infty \dd t\, \mathscr{X}_t\Big)\,\prod_{i=0}^\infty F(\mu_1,\mu_2, (\Delta\mathscr{X}_{t_i})^2)\Big).$$
We now recall that the distribution of $(\mathscr{X}_t)_{t\geq 0}$ is the excursion measure of the $\phi$-CSBP in order
to rewrite this equality in a slightly different form. Suppose that $\sum_{k\in K}\delta_{\omega_k}$ is a Poisson point measure
with intensity $\N_0$. The process $(X_t)_{t\geq 0}$ defined by $X_0=1$ and $X_t=\sum_{k\in K} \mathscr{X}_t(\omega_k)$ if $t>0$ is 
then a $\phi$-CSBP started at $1$. Furthermore, the exponential formula for Poisson measures and the last display immediately give
\begin{equation} 
\label{distriple1}\E\Big[\exp\Big(-\lambda\int_0^\infty \dd t\, X_t\Big)\,\prod_{j=0}^\infty F(\mu_1,\mu_2, (\Delta X_{s_j})^2)\Big]=\exp(-v(\lambda,\mu_1,\mu_2))
\end{equation}
where we have written $(s_j)_{j\geq 0}$ for a measurable enumeration of the jumps of $X$. 

Let $t\geq 0$. Using the Markov property of $X$ at time $t$, the left-hand side of \eqref{distriple1} is also equal to
\begin{equation} 
\label{distriple2}\E\Big[\Big(\exp\Big(-\lambda\int_0^t \dd s\, X_s\Big)\,\prod_{j:s_j\leq t} F(\mu_1,\mu_2, (\Delta X_{s_j})^2)\Big) \exp(-v(\lambda,\mu_1,\mu_2)X_t)\Big].
\end{equation}
To simplify notation, we write $v=v(\lambda,\mu_1,\mu_2)$ in the following calculations, which are very similar to the proof
of Proposition 4.8 in \cite{CLG}. We also set, for every $s\geq 0$,
$$V_s:= \exp\Big(-\lambda\int_0^s \dd u\, X_u\Big)\,\prod_{j:s_j\leq s} F(\mu_1,\mu_2, (\Delta X_{s_j})^2).$$
 From the form of the generator of the $\phi$-CSBP, we have
$$e^{-v X_t}= e^{-v } + M_t + \phi(v)\int_0^t X_s\,e^{- v X_s}\,\dd s,$$
where $(M_s)_{s\geq 0}$ is a martingale, which is bounded on every compact time interval. By using the integration by parts formula as in
\cite[formula (28)]{CLG}, we get
$$e^{-v X_t} V_t= e^{-v } + \int_0^t V_{s-}\,\dd M_s + \phi(v)\int_0^t V_sX_s\,e^{- v X_s}\,\dd s + \int_0^t e^{-v X_s}\,\dd V_s.$$
From \eqref{distriple1} and \eqref{distriple2}, we have $\E[e^{-v X_t} V_t]=e^{-v }$. Hence, taking expectations in the
last display, we obtain
$$\phi(v)\,\E\Big[\int_0^t V_sX_s\,e^{- v X_s}\,\dd s\Big] = - \E\Big[\int_0^t e^{-v X_s}\,\dd V_s\Big].$$
Next we observe that
$$\int_0^t e^{-v X_s}\,\dd V_s= -\lambda \int_0^t V_sX_s\,e^{- v X_s}\,\dd s + \sum_{j:s_j\leq t} e^{-v X_{s_j}}\,V_{s_j-}\, (F(\mu_1,\mu_2,(\Delta X_{s_j})^2)-1),$$
and so we get
$$(\phi(v)-\lambda)\,\E\Big[\int_0^t V_sX_s\,e^{- v X_s}\,\dd s\Big]  = - \E\Big[\sum_{j:s_j\leq t} e^{-v X_{s_j}}\,V_{s_j-}\, (F(\mu_1,\mu_2,(\Delta X_{s_j})^2)-1)\Big].$$
We multiply both sides of this identity by $1/t$ and let $t\da 0$. We have first
$$\lim_{t\da 0} \frac{1}{t}\,\E\Big[\int_0^t V_sX_s\,e^{- v X_s}\,\dd s\Big] = e^{-v}.$$
On the other hand, as a consequence of the classical Lamperti representation of
continuous-state branching processes \cite{Lam0,CLU}, we know that the dual predictable projection of the random measure
$$\sum_{i=0}^\infty \delta_{(s_j,\Delta X_{s_j})}(\dd s,\dd x)$$
is the measure $X_s\,\dd s\,\kappa(\dd x)$, where $\kappa(\dd x)=\sqrt{3/2\pi}\,x^{-5/2}\,\mathbf{1}_{\{x>0\}}\dd x$ is the L\'evy measure 
of the L\'evy process appearing in the Lamperti representation of $X$. This implies that
$$\E\Big[\sum_{j:s_j\leq t} e^{-v X_{s_j}}\,V_{s_j-}\, (F(\mu_1,\mu_2,(\Delta X_{s_j})^2)-1)\Big]
=\E\Big[\int_0^t \dd s\, e^{-v X_s}\,V_sX_s\!\int \kappa(\dd x)\, e^{-v x}(F(\mu_1,\mu_2,x^2)-1)\Big].$$
Consequently,
$$\lim_{t\da 0} \frac{1}{t}\,\E\Big[\sum_{j:s_j\leq t} e^{-v X_{s_j}}\,V_{s_j-}\, (F(\mu_1,\mu_2,(\Delta X_{s_j})^2)-1)\Big]
= e^{-v}\,\int \kappa(\dd x)\, e^{-v x}(F(\mu_1,\mu_2,x^2)-1).$$
Finally, we have obtained
$$\phi(v)-\lambda= - \int \kappa(\dd x)\, e^{-v x}(F(\mu_1,\mu_2,x^2)-1).$$
Using the equality $\phi(v)=\int \kappa(\dd x)(e^{-v x}-1+vx)$, straightforward calculations left to the reader show that
$$\phi(v)+ \int \kappa(\dd x)\, e^{-v x}(F(\mu_1,\mu_2,x^2)-1) = \frac{1}{\sqrt{6}}\,h_{\mu_1,\mu_2}(v),$$
where $h_{\mu_1,\mu_2}$ is as in the statement. This proves that $v=v(\lambda,\mu_1,\mu_2)$ solves $h_{\mu_1,\mu_2}(v)=\sqrt{6}\,\lambda$.
Uniqueness is clear since the function $h_{\mu_1,\mu_2}$ is monotone increasing over $[0,\infty)$. \endproof

\begin{corollary}
\label{disLTduration}
For every $\lambda\geq 0$ and $\mu > 0$, we
have
$$\N_0\Big(1-\exp(-\lambda \ll^0-\mu\sigma)\Big)=\left\{ \begin{array}{ll}
\sqrt{2\mu}\,\cos\Big(\frac{2}{3}\,\arccos\Big( \frac{\sqrt{3}\,\lambda}{2 (2\mu)^{3/4}}\Big)\Big)\quad&\hbox{if } \frac{\sqrt{3}\,\lambda}{2 (2\mu)^{3/4}}\leq 1,\\
\noalign{\smallskip}
\sqrt{2\mu}\,\cosh\Big(\frac{2}{3}\,\mathrm{arcosh}\Big( \frac{\sqrt{3}\,\lambda}{2 (2\mu)^{3/4}}\Big)\Big)\quad&\hbox{if } \frac{\sqrt{3}\,\lambda}{2 (2\mu)^{3/4}}\geq 1.
\end{array}
\right.
$$
\end{corollary}

\proof Set $w(\lambda,\mu)=\N_0(1-\exp(-\lambda \ll^0-\mu\sigma))$. Note that $w(\lambda,\mu)\geq \N_0(1-\exp(-\mu\sigma))=\sqrt{\mu/2}$
by \eqref{Laplace-sigma}.
It follows from Proposition \ref{distriple} applied with $\mu_1=\mu_2=\mu$ that $w(\lambda,\mu)$ is the unique solution of the equation
$$4 w^3- 6\mu w +(2\mu)^{3/2} = \frac{3}{2}\lambda^2$$
in $[\sqrt{\mu/2},\infty)$ (note that the left-hand side is a monotone increasing function of $w$ on $[\sqrt{\mu/2},\infty)$). Set $\wt w(\lambda,\mu)= w(\lambda,\mu)/\sqrt{2\mu}$ and $a= \sqrt{3}\,\lambda/(2(2\mu)^{3/4})$. We immediately get that $\wt w(\lambda,\mu)$ is the unique solution of
$$4 \wt w^3 -3\wt w + 1 = 2a^2$$
in  $[1/2,\infty)$. A simple calculation now shows that
$$\wt w = \left\{
\begin{array}{ll}
\cos(\frac{2}{3}\arccos(a))\quad&\hbox{if }a\leq 1,\\
\noalign{\smallskip}
\cosh(\frac{2}{3}\,\mathrm{arcosh}(a))\quad&\hbox{if }a\geq 1,
\end{array}
\right.
$$
solves the preceding equation. This completes the proof. \endproof

We can also derive an explicit formula for $\N_x(1-\exp(-\lambda \ll^0-\mu\sigma))$, for every $x\in\R$, from Corollary \ref{disLTduration}. Fix $x>0$ for definiteness
and argue under the measure $\N_x$. Write $T_0(\w)=\inf\{t\in[0,\zeta_{(\w)}]:\w(t)=0\}$ for any finite path $\w$ and define
$$\y_0=\int_0^\sigma \dd s\,\mathbf{1}_{\{T_0(W_s)=\infty\}}.$$
Also let $(\omega_i)_{i\in I}$ be the excursions outside $(0,\infty)$ defined as in Section \ref{preli}. Then, we have $\N_x$ a.e.
$$\sigma=\y_0 + \sum_{i\in I} \sigma({\omega_i})\;,\quad \ll^0=\sum_{i\in I}\ll^0(\omega_i)$$
where the second equality follows from the proof of Corollary \ref{superBM}. Using the 
special Markov property (with the fact that $\mathcal{Y}_0$ is $\ee_{(0,\infty)}$-measurable), we get
\begin{equation}
\label{tech00}
\N_x\Big(1-\exp(-\lambda \ll^0-\mu\sigma)\Big)
= \N_x\Big(1-\exp\Big(-\mu\y_0 - \z_0\N_0(1-\exp(-\lambda \ll^0-\mu\sigma))\Big).
\end{equation}
On the other hand, Lemma 4.5 in \cite{CLG} shows that, for every $\mu,\theta>0$ such that $\theta\geq \sqrt{\mu/2}$,
\begin{equation}
\label{tech01}
\N_x(1-\exp(-\mu \y_0-\theta\z_0))=
\sqrt{\frac{\mu}{2}}\Bigg( 3\Bigg(\coth\Bigg((2\mu)^{1/4} x +\coth^{-1}\sqrt{\frac{2}{3} +\frac{1}{3}
\sqrt{\frac{2}{\mu}}\,\theta}\,\Bigg)\Bigg)^2 -2 \Bigg)
\end{equation}
with the convention that the right-hand side equals $\sqrt{\mu/2}$ if $\theta=\sqrt{\mu/2}$. 

Taking $\theta=\N_0(1-\exp(-\lambda \ll^0-\mu\sigma))\geq \sqrt{\mu/2}$ in \eqref{tech01}, using the formula of Corollary \ref{disLTduration}, then yields
a (complicated but explicit) expression for $\N_x(1-\exp(-\lambda \ll^0-\mu\sigma))$. 

\begin{corollary}
\label{dispair}
For every $\mu_1,\mu_2\geq 0$, we have
$$\N_0(1-\exp(-\mu_1\sigma_+-\mu_2\sigma_-))= \frac{\sqrt{2}}{3}\,\frac{\mu_1^{3/2}-\mu_2^{3/2}}{\mu_1-\mu_2},$$
with the convention that the right-hand side equals $\sqrt{\mu_1/2}$ if $\mu_1=\mu_2$. 
The distribution of the pair $(\sigma_+,\sigma_-)$ under $\N_0$ has density 
$$g(s_1,s_2):=\frac{1}{2\sqrt{2\pi}}\,(s_1+s_2)^{-5/2}$$
with respect to Lebesgue measure on $(0,\infty)^2$. In particular, the distribution of $\sigma_+$ (or of $\sigma_-$)
under $\N_0$ has density $(3\sqrt{2\pi})^{-1}\,s^{-3/2}$ on $(0,\infty)$. 

\end{corollary}

The form 
of the density $g(s_1,s_2)$ shows that the conditional distribution of $\sigma_+$ knowing that $\sigma=s$ is uniform
over $[0,s]$. This is a well-known fact, which can be derived from the invariance of the CRT under uniform re-rooting 
(see e.g.~\cite[Section 3.2]{Al2}).

\proof The formula for $\N_0(1-\exp(-\mu_1\sigma_+-\mu_2\sigma_-))$ is obtained by solving the equation $h_{\mu_1,\mu_2}(v)=0$. 
We can then verify that the function $g$ satisfies
$$\int_0^\infty\!\!\int_0^\infty \dd s_1\,\dd s_2 \,g(s_1,s_2)\,(1-e^{-\mu_1s_1-\mu_2s_2})= \frac{\sqrt{2}}{3}\,\frac{\mu_1^{3/2}-\mu_2^{3/2}}{\mu_1-\mu_2},$$
which gives the second assertion. \endproof

We finally give an application to super-Brownian motion in the spirit of Corollary \ref{superBM}.

\begin{corollary}
\label{superBM2}
Let $\mathbf X$ be  a one-dimensional super-Brownian motion with branching mechanism
$\psi(u)=2u^2$, such that $\mathbf{X}_0=\alpha\delta_0$. Set 
$$\mathbf{R}_+=\int_0^\infty \dd t\,\langle \mathbf{X}_t,\mathbf{1}_{[0,\infty)}\rangle\;,\quad \mathbf{R}_-=\int_0^\infty \dd t\,\langle \mathbf{X}_t,\mathbf{1}_{(-\infty,0]}\rangle.$$
Then, for every $\mu_1,\mu_2>0$,
$$\E[\exp(-\mu_1\mathbf{R}_+ - \mu_2\mathbf{R}_-)]= \exp\Big( -\frac{\alpha\sqrt{2}}{3}\,\frac{\mu_1^{3/2}-\mu_2^{3/2}}{\mu_1-\mu_2}\Big).$$
\end{corollary}

Given Corollary \ref{dispair}, the proof of Corollary \ref{superBM2} is an immediate application of the Brownian snake
construction of super-Brownian motion along the lines of the proof of Corollary \ref{superBM}. 

\section{Conditional distributions of the local time at $0$}
\label{seccondi}

We will now use the preceding results to recover 
the conditional distribution of $\ll^0$ given $\sigma$, which was first obtained
by Bousquet-M\'elou and Janson \cite{BMJ} with a very different method.

\begin{theorem}
\label{lawLT}
Let $s>0$. Under the probability measure $\N_0(\cdot\mid \sigma=s)$, the local time $\ll^0$
is distributed as $(2^{3/4}/3)\,s^{3/4}\,T^{-1/2}$, where $T$ is a positive stable
variable with index $2/3$, whose Laplace transform is $\E[\exp(-\lambda T)]=\exp(-\lambda^{2/3})$.
\end{theorem}

\rem In Corollary 3.4 of \cite{BMJ}, the constant $2^{3/4}/3$ is replaced by $2^{1/4}/3$. This is due to
a different normalization: In \cite{BMJ} (as in \cite{Al2}) the random function coding the genealogy of ISE is
twice the Brownian excursion, and it follows that our random variable $\ll^0$ is distributed 
under $\N_0(\cdot\,|\,\sigma=1)$ as $\sqrt{2}$ times the quantity $f_{\mathrm{ISE}}(0)$ considered in \cite{BMJ}.

\smallskip

The occurence of a stable variable with index $2/3$ in Theorem \ref{lawLT} is of course reminiscent of
Corollary \ref{superBM} above. It would be very interesting to establish a direct connection
between this corollary and Theorem \ref{lawLT}. 

\proof From the scaling properties of the end of Section \ref{preli}, it is enough to treat the case $s=1$. Recall the notation $v(\lambda,\mu_1,\mu_2)$ in Proposition \ref{distriple}. For every $\lambda\geq0$, set
$$F(\lambda) := 2\,\N_0\Big(e^{-\sigma/2} (1-e^{-\lambda\ll^0})\Big)= 2\,v(\lambda,\frac{1}{2},\frac{1}{2}) - 1,$$
where the second equality holds because $\N_0(1-\exp(-\sigma/2))=1/2$. The function $F$ is continuous and
vanishes at $0$. As a straightforward consequence of Proposition \ref{distriple}, we have for every $\lambda\geq 0$,
\begin{equation}
\label{equaF}
F(\lambda) = \lambda\,\sqrt{\frac{3}{3+F(\lambda)}}.
\end{equation}
In particular, the right derivative of $F$ at $0$ is $1$, and consequently $\N_0(\ll^0\,\exp(-\sigma/2))=1/2$. The fact that
$\N_0(\ll^0\,\exp(-\sigma/2))$ is finite allows us to make sense of $F(\lambda)$ for every $\lambda\in \C$
such that $\mathrm{Re}(\lambda)\geq 0$, and the restriction of $F$ to $\{\lambda\in\C:\mathrm{Re}(\lambda)>0\}$
is analytic. 

Set $\psi(z)=\sqrt{3/(3+z)}$ so that $\psi$ is analytic on a neighborhood of $0$ in $\C$. Since $\psi(0)\not=0$, 
we can find an analytic function $G$ defined on a neighborhood of $0$ such that $z\psi(G(z))=G(z)$
for $|z|$ small enough. By \eqref{equaF}, we must have $F(z)=G(z)$ for $\mathrm{Re}(z)>0$ and $|z|$ small, and this
means that $F$ can be extended to an analytic function on a neighborhood of $0$.
By the Lagrange inversion theorem, we have then, for every integer $n\geq 1$,
$$[z^n]F(z)=\frac{1}{n}[z^{n-1}]\psi(z)^n= \frac{3^{n/2}}{n!}\,\frac{\dd^{n-1}(3+z)^{-n/2}}{\dd z^{n-1}}\Big|_{\displaystyle z=0}= \frac{(-1)^{n-1}3^{1-n}}{n!}
\frac{\Gamma(\frac{3n}{2}-1)}{\Gamma(\frac{n}{2})},$$
using the standard notation $[z^n]F(z)$ for the coefficient of $z^n$ in the series expansion of $F(z)$ near $0$.
On the other hand, the fact that the function $z\mapsto F(z)$ is analytic in a neighborhood of $0$ implies that all moments
$\N_0((\ll^0)^ne^{-\sigma/2})$, $n\geq 1$, are finite and given by
\begin{equation}
\label{moment-local}
\N_0((\ll^0)^ne^{-\sigma/2})=\frac{1}{2}\,(-1)^{n-1}n!\,\times [z^n]F(z)= \frac{1}{2}\,3^{1-n}\,\frac{\Gamma(\frac{3n}{2}-1)}{\Gamma(\frac{n}{2})}.
\end{equation}
To complete the proof, we use a scaling argument. We recall that the distribution of $\ll^0$ under $\N_0(\cdot\mid\sigma=s)$
coincides with the distribution of $s^{3/4}\ll^0$ under $\N_0(\cdot\mid\sigma=1)$. It follows that
$$\N_0((\ll^0)^ne^{-\sigma/2})=\int_0^\infty \frac{\dd s}{2\sqrt{2\pi s^3}}\, e^{-s/2} \N_0\Big(s^{3n/4}(\ll^0)^n\,\Big|\,\sigma=1\Big)
= \frac{2^{\frac{3n}{4}-2}}{\sqrt{\pi}}\, \Gamma(\frac{3n}{4}-\frac{1}{2})\times \N_0\Big((\ll^0)^n\,\Big|\,\sigma=1\Big).$$

By combining the last two displays and using the duplication formula for the Gamma function, we arrive at
$$\N_0\Big((\ll^0)^n\,\Big|\,\sigma=1\Big) = \frac{\sqrt{\pi}\,3^{1-n}}{2^{\frac{3n}{4}-1}}\,\frac{\Gamma(\frac{3n}{2}-1)}{\Gamma(\frac{n}{2})\Gamma(\frac{3n}{4}-\frac{1}{2})}
= \frac{2^\frac{3n}{4}}{3^n}\,\frac{\Gamma(\frac{3n}{4}+1)}{\Gamma(\frac{n}{2}+1)}= \Big( \frac{2^{3/4}}{3}\Big)^n\E[T^{-n/2}],$$
where $T$ is as in the theorem (to check the last equality, write $T^{-n/2}= (\Gamma(n/2))^{-1}\!\int_0^\infty \dd s \,s^{n/2-1}e^{-sT}$). 
The growth of the moments of the distribution of $T^{-1/2}$ ensures that this distribution is characterized by its moments, which
completes the proof. \endproof

\rem Rather than using the Lagrange inversion theorem, we could have derived formula \eqref{moment-local} for the
moments $\N_0((\ll^0)^ne^{-\sigma/2})$ from a series expansion of the quantity $\N_0(1-\exp(-\lambda\ll^0-\sigma/2))$
as given in Corollary \ref{disLTduration}. This would still have required some calculations. We preferred to
use the previous method because it also serves as a prototype for the proof of the (more delicate) Theorem \ref{lawLT+} below. 

\medskip

Proposition \ref{distriple} can also be used to derive the conditional distribution of $\ll^0$
given $\sigma_+$. Perhaps surprisingly, this distribution turns out again to be remarkably
simple.

\begin{theorem}
\label{lawLT+}
Let $s>0$. Under the probability measure $\N_0(\cdot\mid \sigma_+=s)$, the local time $\ll^0$
is distributed as $(2^{9/4}/3)\,s^{3/4}\,D\,T^{-1/2}$, where 
the random variables $D$ and $T$ are independent, $T$ is a positive stable
variable with index $2/3$, whose Laplace transform is $\E[\exp(-\lambda T)]=\exp(-\lambda^{2/3})$,
and $D$ has density $2x\,\mathbf{1}_{[0,1]}(x)$ with respect to Lebesgue measure on $\R_+$.
\end{theorem}

\proof It is enough to treat the case $s=1$. For every $\lambda\geq 0$,
set
$$F_+(\lambda) := \N_0\Big(1-\exp(-\lambda\ll^0-\frac{1}{2}\sigma_+)\Big)= v(\lambda,\frac{1}{2},0)$$
with the notation of Proposition \ref{distriple}. As for Theorem \ref{lawLT}, the strategy of the proof 
is to compute the coefficient $[\lambda^n]F_+(\lambda)$ in two different ways. Unfortunately, the 
details of the argument are more involved than in the proof of Theorem \ref{lawLT}. 

By Proposition \ref{distriple}, we have
$$(2 F_+(\lambda)-1)\sqrt{F_+(\lambda)+1} + 2\,F_+(\lambda)^{3/2} =\sqrt{6}\,\lambda.$$
We cannot apply directly the Lagrange inversion theorem, but the idea will be to find
a rational parametrization of the preceding equation (see e.g. \cite[Section 3]{BM}). It follows from
the last display that we have $P(F_+(\lambda),\lambda)=0$, where
$$P(y,z)= 96\,y^3z^2 - 36\,z^4 -36\,yz^2+12 \,z^2-9\,y^2+6\,y-1, \qquad y,z\in\C.$$
We now introduce\footnote{The functions $Q$ and $\psi$ have been found using the Maple package {\it algcurve}} the rational functions
$$Q(z)=- \frac{1}{124416}\,z^3 + \frac{1}{48}\,z,\qquad R(z)=\frac{1}{3456}\,z^2-\frac{1}{2} +\frac{216}{z^2},$$
which satisfy $P(R(z),Q(z))=0$ for every $z\in\C\backslash\{0\}$. We have $Q^{-1}(0)=\{-36\sqrt{2},0,36\sqrt{2}\}$, and 
the derivative $Q'$ does not vanish on $Q^{-1}(0)$. It follows that we can find $r_0>0$ and
three analytic functions $\gamma_1,\gamma_2,\gamma_3$ defined on the disk $\D_{r_0}=\{z\in\C:|z|<r_0\}$ and with
disjoint ranges, such that $\gamma_1(0)=-36\sqrt{2},\,\gamma_2(0)=0$, $\gamma_3(0)=36\sqrt{2}$ and for every
$z\in \D_{r_0}$, $Q^{-1}(z)=\{\gamma_1(z),\gamma_2(z),\gamma_3(z)\}$. Note that $R(\gamma_1(0))=1/3=R(\gamma_3(0))$
and $R'(\gamma_1(0))=-\sqrt{2}/54=-R'(\gamma_3(0))$. Also the fact that $Q(\gamma_i(z))=z$ readily implies that
$\gamma'_1(0)=\gamma'_2(0)=-24$.

Since $P(R(z),Q(z))=0$ for every $z\in\C\backslash\{0\}$, we get that $P(R(\gamma_i(z)),z)=0$
for every $i\in\{1,2,3\}$ and $z\in \D_{r_0}\backslash\{0\}$. We claim that $F_+(\lambda)=R(\gamma_1(\lambda))$
for $\lambda>0$ small enough. To see this, observe that for $z\not=0$ and $|z|$ small enough, then the quantities
$R(\gamma_i(z))$, $i\in\{1,2,3\}$, are distinct. Indeed, since $|R(y)|\la \infty$ as $|y|\to 0$ it is
clear that $R(\gamma_2(z))$ is distinct from $R(\gamma_1(z))$ and $R(\gamma_3(z))$
when $|z|$ is small, and on the other hand, the properties $\gamma'_1(0)=\gamma'_3(0)\not =0$
and $R'(\gamma_1(0))=-R'(\gamma_3(0))\not =0$ imply that $R(\gamma_1(z))\not =R(\gamma_3(z))$
when $|z|$ is small. Hence, for $z\not=0$ and $|z|$ small enough, the numbers $R(\gamma_i(z))$, $i\in\{1,2,3\}$, are 
three distinct roots of  $P(y,z)$ viewed as 
a polynomial of degree $3$ in $y$. Since we know that $P(F_+(\lambda),\lambda)=0$, it follows that
$F_+(\lambda)\in\{R(\gamma_1(\lambda)),R(\gamma_2(\lambda)),R(\gamma_3(\lambda))\}$ for $\lambda>0$ small. 
The case $F_+(\lambda)=R(\gamma_2(\lambda))$ is clearly excluded for $\lambda$ small, and since 
$F_+(\lambda)$ is a monotone increasing function of $\lambda$, noting that $\gamma'_1(0)R'(\gamma_1(0))>0$
whereas $\gamma'_3(0)R'(\gamma_3(0))<0$,
we get our claim $F_+(\lambda)=R(\gamma_1(\lambda))$ for $\lambda>0$ small.

In particular, we can extend $F_+$ to an analytic function in the neighborhood of $0$, and we will then use the Lagrange
inversion theorem to determine the coefficients of the Taylor expansion of $F_+$. To simplify notation,
we set $\wt F_+(\lambda)=F_+(\lambda)-1/3$, $\wt \gamma(z)=\gamma_1(z)+36\sqrt{2}$
and for every $\lambda\geq 0$,
$$\wt R(\lambda)=R(\lambda-36\sqrt{2}) - 1/3.$$
Then, for $\lambda>0$ small, we have
\begin{equation}
\label{Lagra1}
\wt F_+(\lambda)=F_+(\lambda)-\frac{1}{3}=R(\gamma_1(\lambda))-\frac{1}{3}=\wt R(\wt\gamma(\lambda)).
\end{equation}
On the other hand, the property $Q(\gamma_1(z))=z$ for $|z|<r_0$ shows that
\begin{equation}
\label{Lagra2}
\wt\gamma(\lambda)=\lambda\,\wt \psi(\wt\gamma(\lambda)),
\end{equation}
with 
$$\wt\psi(\lambda)= - \frac{124416}{(36\sqrt{2}-\lambda)(72\sqrt{2}-\lambda)}.$$
By \eqref{Lagra1}, \eqref{Lagra2} and the Lagrange inversion theorem, we get for every $n\geq 1$,
$$[\lambda ^n]F_+(\lambda )=[\lambda ^n]\wt F_+(\lambda )=\frac{1}{n}[\lambda ^{n-1}](\wt R'(\lambda )\wt\psi(\lambda )^n).$$
Note that
\begin{align*}
&\wt R'(72 \sqrt{2}\lambda)= R'(72 \sqrt{2}(\lambda-1))= -\frac{\sqrt{2}}{48}(1-2\lambda) +\frac{1}{216\sqrt{2}}(1-2\lambda)^{-3}\\
&\wt\psi(72 \sqrt{2}\lambda)= -\frac{24}{(1-\lambda)(1-2\lambda)},
\end{align*}
from which it follows that
\begin{align*}&[\lambda^{n-1}](\wt R'(72\sqrt{2}\lambda)\wt\psi(72\sqrt{2}\lambda)^n)\\
&\quad=(-24)^n\times\Bigg( \Big(\frac{-\sqrt{2}}{48} [\lambda^{n-1}]\Big((1-2\lambda)^{-n+1}(1-\lambda)^{-n}\Big)
+\frac{1}{216\sqrt{2}}[\lambda^{n-1}]\Big((1-2\lambda)^{-n-3}(1-\lambda)^{-n}\Big)\Bigg),
\end{align*}
and finally
\begin{equation}
\label{lawLT+0}
[\lambda^n]F_+(\lambda)= \frac{(-1)^n}{n}\,(3\sqrt{2})^{-n}\,\Big(-3[\lambda^{n-1}]\Big((1-2\lambda)^{-n+1}(1-\lambda)^{-n}\Big)
+\frac{1}{3} [\lambda^{n-1}]\Big((1-2\lambda)^{-n-3}(1-\lambda)^{-n}\Big)\Big).
\end{equation}
To compute the right-hand side, we observe that, for every integers $m\geq 0$, $k\geq 1$ and $\ell\geq 1$,
we have
$$[\lambda^m](1-2\lambda)^{-k}(1-\lambda)^{-\ell} = 2^m\,{m+k-1\choose m}\,_2F_1(-m,\ell;-m-k+1;\frac{1}{2}),$$
where $_2F_1$ stands for the Gauss hypergeometric function. This equality is easily checked by a direct
calculation, noting that the hypergeometric series reduces to a finite sum in the case we are considering
It follows that, for every $n\geq 2$,
\begin{align}
\label{lawLT+1}
&[\lambda^{n-1}]\Big((1-2\lambda)^{-n+1}(1-\lambda)^{-n}\Big)= 2^{n-1} {2n-3\choose n-1}\,
_2F_1(-n+1,n;-2n+3;\frac{1}{2})\\
\label{lawLT+2}
&[\lambda^{n-1}]\Big((1-2\lambda)^{-n-3}(1-\lambda)^{-n}\Big)= 2^{n-1} {2n+1\choose n-1}\,
_2F_1(-n+1,n;-2n-1;\frac{1}{2}).
\end{align}
Fortunately, Bailey's theorem (see \cite[Theorem 3.5.4 (ii)]{AAR}) gives an explicit formula for 
$_2F_1(a,1-a;b;\frac{1}{2})$ in terms of a ratio of products of values of the
Gamma function, which we can apply here. Using also Euler's reflection formula
$\Gamma(z)\Gamma(1-z)=\pi/\sin(\pi z)$ to eliminate the poles of the Gamma function, 
we arrive at
\begin{align*}
_2F_1(-n+1,n;-2n+3;\frac{1}{2})&=\frac{\Gamma(\frac{n}{2}-\frac{1}{2})\,\Gamma(\frac{3n}{2}-1)}
{\Gamma(n-\frac{1}{2})\,\Gamma(n-1)}= \frac{2^{2n-3}}{\sqrt{\pi}}\,\frac{\Gamma(\frac{n}{2}-\frac{1}{2})\,\Gamma(\frac{3n}{2}-1)}
{\Gamma(2n-2)}\\
_2F_1(-n+1,n;-2n-1;\frac{1}{2})&=\frac{\Gamma(\frac{n}{2}+\frac{3}{2})\,\Gamma(\frac{3n}{2}+1)}
{\Gamma(n+\frac{3}{2})\,\Gamma(n+1)}=\frac{2^{2n+1}}{\sqrt{\pi}}\,\frac{\Gamma(\frac{n}{2}+\frac{3}{2})\,\Gamma(\frac{3n}{2}+1)}
{\Gamma(2n+2)},
\end{align*}
where we applied the duplication formula for the Gamma function, and we recall that we assume $n\geq 2$.
Using \eqref{lawLT+1} and \eqref{lawLT+2}, we get from \eqref{lawLT+0} that
\begin{align*}
[\lambda^n]F_+(\lambda)&= \frac{(-1)^{n+1}}{n!}\,(3\sqrt{2})^{-n}\,\frac{2^{3n}}{\sqrt{\pi}}\,\Big(\frac{3}{16}\, \frac{\Gamma(\frac{n}{2}-\frac{1}{2})\,\Gamma(\frac{3n}{2}-1)}
{\Gamma(n-1)} -\frac{1}{3}\,\frac{\Gamma(\frac{n}{2}+\frac{3}{2})\,\Gamma(\frac{3n}{2}+1)}{\Gamma(n+3)}\Big)\\
&=\frac{(-1)^{n+1}}{n!}\,(3\sqrt{2})^{-n}\,\frac{2^{3n}}{\sqrt{\pi}}\,\frac{\Gamma(\frac{n}{2}+\frac{1}{2})\,\Gamma(\frac{3n}{2}-1)}{\Gamma(n)}
\Big( \frac{3}{8} - \frac{1}{3}\times\frac{(\frac{3n}{2}-1)(\frac{n}{2} +\frac{1}{2})\frac{3n}{2}}{(n+2)(n+1)n}\Big)\\
&=\frac{(-1)^{n+1}}{n!}\,(3\sqrt{2})^{-n}\,\frac{2^{3n}}{\sqrt{\pi}}\,\frac{1}{n+2}\,\frac{\Gamma(\frac{n}{2}+\frac{1}{2})\,\Gamma(\frac{3n}{2}-1)}{\Gamma(n)}\\
&=\frac{(-1)^{n+1}}{n!}\,(3\sqrt{2})^{-n}\,2^{2n+1}\,\frac{1}{n+2}\,\frac{\Gamma(\frac{3n}{2}-1)}{\Gamma(\frac{n}{2})}.
\end{align*}
We have assumed $n\geq 2$, but a direct calculation from \eqref{lawLT+0} shows that the last line of the preceding display also gives
the correct value $[\lambda]F_+(\lambda)=4\sqrt{2}/9$ for $n=1$. Similarly as in the proof of 
Theorem \ref{lawLT}, we conclude that, for every $n\geq 1$,
$$\N_0\Big( (\ll^0)^n\,e^{-\sigma_+/2}\Big)= \Big(\frac{2\sqrt{2}}{3}\Big)^n\,\frac{2}{n+2}\,\frac{\Gamma(\frac{3n}{2}-1)}{\Gamma(\frac{n}{2})}.$$
On the other hand, the same scaling argument as in the proof of Theorem \ref{lawLT} (using now the fact that the density
of $\sigma_+$ under $\N_0$ is $(3\sqrt{2\pi})^{-1}s^{-3/2}$) gives
$$\N_0((\ll^0)^ne^{-\sigma_+/2})=\int_0^\infty \frac{\dd s}{3\sqrt{2\pi s^3}}\, e^{-s/2} \N_0\Big(s^{3n/4}(\ll^0)^n\,\Big|\,\sigma_+=1\Big)
= \frac{2^{\frac{3n}{4}-1}}{3\sqrt{\pi}}\, \Gamma(\frac{3n}{4}-\frac{1}{2})\, \N_0\Big((\ll^0)^n\,\Big|\,\sigma_+=1\Big).$$
It follows that
$$\N_0\Big((\ll^0)^n\,\Big|\,\sigma_+=1\Big)=3\sqrt{\pi}\Big(\frac{2\sqrt{2}}{3}\Big)^n2^{-\frac{3n}{4}+1}\,\frac{2}{n+2}\,
\frac{\Gamma(\frac{3n}{2}-1)}{\Gamma(\frac{n}{2})\Gamma(\frac{3n}{4}-\frac{1}{2})}= \Big(\frac{2^{9/4}}{3}\Big)^n \, \frac{2}{n+2}\, 
\frac{\Gamma(\frac{3n}{4}+1)}{\Gamma(\frac{n}{2}+1)}.
$$
The right-hand side is the $n$-th moment of $(2^{9/4}/3)\,D\,T^{-1/2}$, where the pair $(D,T)$ is as in the theorem.
This completes the proof. \endproof

\noindent{\bf Interpretation in random geometry.} We now explain briefly how both theorems of this section can be interpreted 
in the setting of continuous models of random geometry. It is best 
to start with the discrete picture of planar quadrangulations. For every integer $n\geq 1$, let $Q_n$ be a uniformly distributed 
rooted and pointed quadrangulation with $n$ faces. The fact that $Q_n$ is pointed means that (in addition to the root edge) there is a distinguished vertex denoted by $\partial$.  Write $\dg$ for the graph distance on 
the vertex set $V(Q_n)$ of $Q_n$. The Schaeffer bijection (see e.g.~\cite[Section 5]{LGM})
allows us to code $Q_n$ by a uniformly distributed labeled tree with $n$ edges, which we denote by $T_n$, and a sign $\epsilon_n\in\{-1,1\}$. Here a labeled tree is a (rooted) plane tree whose vertices are assigned integer labels $\ell_v$, in such a way that the label of the 
root vertex $\rho$ of the tree is $\ell_\rho=0$ and the labels of two adjacent vertices differ by at most $1$ in absolute value. Furthermore 
the set $V(Q_n)\backslash \{\partial\}$ is canonically identified with $V(T_n)$, where $V(T_n)$ denotes the vertex set of $T_n$. Through this identification, the graph distance $\dg(\partial,v)$ between $\partial$ and another 
vertex $v$ of $Q_n$ can be expressed as $\ell_v-\min\{\ell_w:w\in V(T_n)\}+1$. Now consider the set $S_n=\{v\in V(Q_n):\dg(\partial,v)=\dg(\partial,\rho)\}$
of all vertices $v$ of $Q_n$ that are at the same distance as $\rho$ from the distinguished vertex $\partial$ 
(here we view $\rho$ as a vertex of $Q_n$ thanks to the preceding identification). From the previous observations,
$S_n$ is identified to $\{v\in V(T_n): \ell_v=0\}$. It then follows 
from \cite[Theorem 3.6]{BMJ} that the distribution of $n^{-3/4}\# S_n$ converges as $n\to \infty$ to the distribution of
$2^{-1/4}3^{-1/2}\ll^0$ under $\N_0(\cdot\,|\,\sigma=1)$, which is given in Theorem \ref{lawLT}. 

Consider then the (standard) Brownian map $(\bm_\infty,D)$. This is a random compact metric space that can be constructed from Brownian motion indexed by the Brownian tree, which we denote here by $(V_a)_{a\in \t_\zeta}$ as in Section \ref{intro} above,
under the probability measure $\N_0(\cdot\,|\,\sigma=1)$ --- see e.g. the introduction of \cite{Uniqueness} for details. In this construction, the space $\bm_\infty$ is
obtained as a quotient space of $\t_\zeta$, and comes
with two distinguished points, namely the point $\rho$ corresponding to the root of $\t_\zeta$, and 
another point denoted by $x_*$ in \cite{Uniqueness}, which corresponds to the point of $\t_\zeta$ where $V_a$ achieves its minimum.
Note that $\rho$ and $x_*$ can be viewed as independently and uniformly distributed on $\bm_\infty$.
The ``sphere'' $\{x\in\bm_\infty:D(x_*,x)=D(x_*,\rho)\}$ then corresponds to $\{a\in \t_\zeta:V_a=0\}$, and so the local time $\ll^0$
is naturally interpreted as the ``measure'' of this sphere (here the word measure should refer to a suitable Hausdorff measure, although this has not been justified rigorously). This interpretation is made very plausible by the 
discrete result for quadrangulations described above.

To get a similar interpretation for Theorem \ref{lawLT+}, we consider the free Brownian map $(M,\Delta)$, which is 
the scaling limit of quadrangulations distributed according to Boltzmann weights and can 
again be constructed from Brownian motion indexed by the Brownian tree, but now under the $\sigma$-finite measure $\N_0$ (see e.g.~\cite[Section 3]{Disks}).
As in the case of the standard Brownian map, the space $M$ is defined as a quotient space of $\t_\zeta$ and
comes with two distinguished points 
denoted by $\rho$ and $x_*$. Furthermore, the sphere $\{x\in M:\Delta(x_*,x)=\Delta(x_*,\rho)\}$ corresponds to $\{a\in \t_\zeta:V_a=0\}$,
and the ball $\{x\in M:\Delta(x_*,x)\leq \Delta(x_*,\rho)\}$ corresponds to $\{a\in \t_\zeta:V_a\leq 0\}$. So Theorem \ref{lawLT+}
can be viewed as providing the conditional distribution of the measure of the sphere $\{x\in M:\Delta(x_*,x)=\Delta(x_*,\rho)\}$ given the volume of the ball 
it encloses.

\medskip

\noindent{\bf Acknowledgement.} We thank Nicolas Curien for suggesting the use of the Maple package {\it algcurve}
in the proof of Theorem \ref{lawLT+}.

\smallskip
\noindent{\bf Note added in proof.} Alin Bostan (personal communication) has pointed out the following simplification in the proof
of Theorem \ref{lawLT+}. Since the function $F_+$ is algebraic, a classical theorem says
that it is also $D$-finite, hence satisfies a linear differential equation whose coefficients
are polynomials in $z$. This differential equation can be written explicitly and yields
a recursive equation for the coefficients $[\lambda^n]F_+(\lambda)$, which then leads
to our explicit formula for these coefficients.

\end{document}